\documentclass[11pt,a4paper,leqno,fleqn]{article}
\usepackage{graphicx}
\usepackage{ascmac}
\usepackage{hyperref}
\usepackage{amsmath,amssymb}
\usepackage{amsthm}
\usepackage{newpxtext,newpxmath}
\title {Research for a set of particular\\ primitive Pythagorean triples }
\newcommand{\ebrak}[1]{\left( \begin{array}{@{}c@{}} {#1} \end{array} \right)}
\author {Yasushi Ieno}
\date{}
\begin{document}
\maketitle 
\noindent Abstract. The set of terms of an infinite sequence expressed by a recurrence relation is equal to the set of maximum numbers of all primitive Pythagorean triples such that the difference between the two non-maximum numbers is 1, which Cimmino showed. By reference to Cimmino's researches, we show a recurrence relation for the case where the difference between the two non-maximum numbers is not 1, but 7.\\

\noindent Key Words and Phrases. recurrence relation, primitive Pythagorean triple
\\\\
0. Introduction\\

The set of terms of the infinite sequence $\rm \{a_n\}_{n=1,2,\cdots{,}\infty}$ expressed by the recurrence formula $\rm a_{n+1}{=}6a_n {-}a_{n-1}$, where n$\geq$1, starting with terms $\rm a_0{=}1$ and $\rm a_1{=}5$, is equal to the set of the maximum number z's of primitive Pythagorean triples $\rm \{ (x,y,z)| x^2{+}y^2{=}z^2; y{=}x{+}1; x,y,z{\in}N^{+}; \\gcd(x,y,z)=1\}$, which was shown by Cimmino [1]. 

Cimmino [1] first showed that all the terms of the sequence represented by the above-mentioned recurrence relation are equal to z of some element of the above-mentioned set of primitive Pythagorean triples. 

Next, by using parameters effectively, it was shown that the set consisting of each maximum number of any element of this primitive Pythagorean triple set has a one-to-one correspondence with the set consisting of the terms of the sequence represented by the recurrence relation.

By reference to Cimmino's researches [1] we will consider the case where the difference between the two non-maximum numbers is 7, not 1. Then we are starting by predicting the recurrence relation.
\\\\
1. Predicting the recurrence relation.\\

In this paper we consider a set of primitive Pythagorean triples as follows,
\begin{eqnarray}
\lefteqn {\rm \{ (x,y,z)| x^2{+}y^2{=}z^2; y{=}x{+}7; x,y,z{\in}N^{+}; gcd(x,y,z){=}1\} } 
\end{eqnarray}

In [1], Cimmino locally extended the bounds of x to 0 for the case y=x+1. 

Then (x, y, z)=(0, 1, 1) is a 'virtual' primitive Pythagorean triple, holding $\rm x^2{+}y^2{=}z^2$, as $\rm a_0$=1, one of the initial values of the recurrence formula for the case y=x+1. \\

Similarly we allow x to be negative for the case y=x+7. 

For now (x, y, z)=(-4, 3, 5) is a 'virtual' primitive Pythagorean triple, clearly holding $\rm y{=}x{+}7$ and $\rm x^2{+}y^2{=}z^2$. 

Then $\rm a_0$=5, one of the initial values of the recurrence relation.

By a short calculation, we see that (x,y,z)=(5,12,13) is the smallest real primitive Pythagorean triple for the case y=x+7. So suppose $\rm a_1$=13.

Then we search for the recurrence relation for now as follows.
\begin{eqnarray} 
\lefteqn {\rm a_{n{+}1}{=}Aa_n{+}Ba_{n{-}1} }\\
\lefteqn {\rm \ \ \ where\ A\ and\ B\ are\ constant\ values\ and\ n{\geq}1,} \nonumber \\ 
\lefteqn {\rm \ \ \ starting\ with\ terms\ a_0{=}5\ and\ a_1{=}13.} \nonumber
\end{eqnarray}
\qquad
Now we denote the solutions of a quadratic equation $\rm x^{2}{=}Ax{+}B$ by $\rm x_{\pm}$.

Then $\rm {x_{+}}{=}{\dfrac{A{+}\sqrt{A^2{+}4B}}{2}}$, $\rm {x_{-}}={\dfrac{A{-}\sqrt{A^2{+}4B}}{2}}$ and $\rm a_n$ can be written as below.
\begin{eqnarray}
\lefteqn {\rm a_n{=}a{x_{+}}^n{+}b{x_{-}}^n} \ \ \ \ \ \ \ \ \ \ \ \ \ \ \ \ \ \ \ \ \ \ \ \ \ \ \ \ \ \ \ \ \ \ \ \ \ \ \ \ 
\end{eqnarray}

Substituting n=0 and n=1 into (3), we obtain a linear simultaneous equation for a and b.\\
$$ \left\{
\begin{aligned}
\rm a & \ \rm {+} & \rm b & = & \rm a_0 = 5\ \ \ \ \\
\rm a{x_{+}} & \ \rm {+} & \rm b{{x_{-}}} & = & \rm a_1 = 13\ \ 
\end{aligned}
\right.
$$
\qquad 
Solving this equation, 

$\rm \ \ \ \ \ \ \ a={\dfrac{26{-}5A{+}5\sqrt{A^2{+}4B}}{2\sqrt{A^2{+}4B}}}$, $\rm b={\dfrac{5A{-}26{+}5\sqrt{A^2{+}4B}}{2\sqrt{A^2{+}4B}}}$
\\

If there exists an element (x,y,z) of the set expressed at (1) where $\rm z{=}a_n$, \\
then there exists some positive integer m where\\ 
\begin{eqnarray}
\lefteqn {\rm \ \ \ \ m^2{+}(m{+}7)^2{=}a_n^2}  
\end{eqnarray}
\qquad
Transforming (4), 

$\rm {\ \ \ \ \ \ \ \dfrac{a_n^2{-}49}{2}{=}m(m{+}7)}$

$\rm {\ \ \ \ \ \ \  \ebrak{\dfrac{\sqrt{\rm a_n^2{-}49}{-}7}{2}}\ebrak{\dfrac{\sqrt{\rm a_n^2{-}49}{-}7}{2}{+}7}{=}m(m{+}7)}$

So there exists a positive integer $\rm m_n$ which holds

$\rm {\ \ \ \ \ \ \ m_n{=}\dfrac{\sqrt{2a_n^2{-}49}{-}7}{2}}$

Therefore it is followed by
\begin{eqnarray}
\lefteqn {\rm {\ \ \ \ \ \ {2a_n^2{-}49}{=}{d_n}^2}} \\
\lefteqn {\rm {\ \ \ \ \ \ \ \ \ \ \ ,\ where\ {d_n}\ is \ a \ positive\ odd\ integer\ greater\ than\ 7.}} \nonumber
\end{eqnarray}

According to (3), \\
$\rm a_n{=}\ebrak{\dfrac{\rm 5\sqrt{A^2{+}4B}{-}(5A{-}26)}{\rm 2\sqrt{A^2{+}4B}}}\ebrak{\dfrac{\rm A{+}\sqrt{A^2{+}4B}}{\rm 2}}^n{+}\ebrak{\dfrac{\rm 5\sqrt{A^2{+}4B}{+}(5A{-}26)}{\rm 2\sqrt{A^2{+}4B}}}\ebrak{\dfrac{\rm A{-}\sqrt{A^2{+}4B}}{\rm 2}}^n$ 
\\
$\rm 2{a_n}^2{-}49{=}\ebrak{\ebrak{\dfrac{\rm 5\sqrt{A^2{+}4B}{-}(5A{-}26)}{\rm \sqrt{2}\sqrt{A^2{+}4B}}}\ebrak{\dfrac{\rm A{+}\sqrt{A^2{+}4B}}{\rm 2}}^{\rm n}{+}\ebrak{\dfrac{\rm 5\sqrt{A^2{+}4B}{+}(5A{-}26)}{\rm \sqrt{2}\sqrt{A^2{+}4B}}}\ebrak{\dfrac{\rm A{-}\sqrt{A^2{+}4B}}{\rm 2}}^{\rm n}}^2{-}49$

It would be convenient if $\rm{2a_n^2{-}49}$ can be described in the squared form of an expression for any n.

So we assume $\rm{2a_n^2{-}49}$ as follows.\\
$\rm {d_n}^2{=}2{a_n}^2{-}49{=}\ebrak{\ebrak{\dfrac{\rm 5\sqrt{A^2{+}4B}{-}(5A{-}26)}{\rm \sqrt{2}\sqrt{A^2{+}4B}}}\ebrak{\dfrac{\rm A{+}\sqrt{A^2{+}4B}}{\rm 2}}^{\rm n}{-}\ebrak{\dfrac{\rm 5\sqrt{A^2{+}4B}{+}(5A{-}26)}{\rm \sqrt{2}\sqrt{A^2{+}4B}}}\ebrak{\dfrac{\rm A{-}\sqrt{A^2{+}4B}}{\rm 2}}^{\rm n}}^2$
\\
$\rm {\ \ \ \ \ \ \ \ \ \ \ \ \ \ \ \ 4\ebrak{\dfrac{\rm 5\sqrt{A^2{+}4B}{-}(5A{-}26)}{\rm 2\sqrt{A^2{+}4B}}}\ebrak{\dfrac{\rm A{+}\sqrt{A^2{+}4B}}{\rm 2}}^n\ebrak{\dfrac{\rm 5\sqrt{A^2{+}4B}{+}(5A{-}26)}{\rm 2\sqrt{A^2{+}4B}}}\ebrak{\dfrac{\rm A{-}\sqrt{A^2{+}4B}}{\rm 2}}^n}$ 
\\
$\rm {\ \ \ \ \ \ \ \ \ \ \ \ \ \ {=}2\ebrak{\dfrac{\rm 25{A^2{+}4B}{-}(5A{-}26)^2}{\rm {A^2{+}4B}}}\ebrak{\rm {-}B}^n{=}49}$

For constancy of the left side regardless of n, B must be equal to ${-}1$, then

$\rm {\ \ \ \ \ \ \ \ \ \ \ \ \ \ 2\ebrak{\dfrac{\rm 25{A^2{-}4}{-}(5A{-}26)^2}{\rm {A^2{-}4}}}{=}49}$\\

$\rm {\ \ \ \ \ \ \ \ \ \ \ \ \ \ \ \ \ \ \ A{=}6,\dfrac{226}{49}}$

If A=$\dfrac{226}{49}$ then

$\ \ \ \ \ \ \ \ \rm a{=}\dfrac{5\sqrt{2}{+}1}{2\sqrt{2}}, b{=}\dfrac{5\sqrt{2}{-}1}{2\sqrt{2}}, x_{+}{=}\dfrac{113{+}72\sqrt{2}}{49}, x_{-}{=}\dfrac{113{-}72\sqrt{2}}{49}.$\\

So $\rm a_2{=}ax_{+}^2{+}bx_{-}^2$

$\ \ \ \ \ \ \ \ \ \rm {=}\ebrak{\dfrac{5\sqrt{2}{+}1}{2\sqrt{2}}}\ebrak{\dfrac{113{+}72\sqrt{2}}{49}}^2{+}\ebrak{\dfrac{5\sqrt{2}{-}1}{2\sqrt{2}}}\ebrak{\dfrac{113{-}72\sqrt{2}}{49}}^2$

$\ \ \ \ \ \ \ \ \ \rm {=}{\dfrac{2693}{49}}\ {\notin}{\rm N^{+}}$; A=$\dfrac{226}{49}$ is invalid.\\

If A=6 then
\begin{eqnarray}
\lefteqn {\rm \ \ \ \ a_{n{+}1}{=}6a_n{-}a_{n{-}1}.} 
\end{eqnarray}
\qquad Now evidently any $\rm a_n$ is a positive integer. And we easily see that $\rm a_2{=}6a_1{-}a_0{=}6{\times}13{-}5{=}73$, and that $73^2{=}48^2{+}55^2{=}48^2{+}(48{+}7)^2$, which results in that (48,55,73) is a primitive Pythagorean triple and an element of the set shown at (1). 73 is a suitable integer for us now.

Besides, 
\begin{eqnarray}
\lefteqn {\rm  {\ {d_n}^2{=}2{a_n}^2{-}49}{=}\ebrak{\ebrak{\dfrac{5\sqrt2{-}1}{2}}(3{+}2\sqrt2)^n{+}\ebrak{\dfrac{5\sqrt2{+}1}{2}}(3{-}2\sqrt2)^n}^2{-}49}\nonumber \\ 
\lefteqn {\rm  \ \ \ \ \ \ \ {=}\ebrak{\ebrak{\dfrac{5\sqrt2{-}1}{2}}(3{+}2\sqrt2)^n{-}\ebrak{\dfrac{5\sqrt2{+}1}{2}}(3{-}2\sqrt2)^n}^2}
\end{eqnarray}

If n${\geq}$1 then $\rm{d_n}^2$ is the square of a positive integer greater than 7.\\

But the current $\rm \{a_n\}_{n=1,2,\cdots{,}\infty}$ is not the required answer.

Now we set focus on a positive integer 17 that is not a term of $\rm \{a_n\}_{n=1,2,\cdots{,}\infty}$, for this sequence monotonically increases with respect to n.

Moreover we easily see that $17^2{=}8^2{+}15^2{=}8^2{+}(8{+}7)^2$, which results in that (8,15,17) is also a primitive Pythagorean triple and an element of the set shown at (1).

Now not only (x, y, z)=(-4, 3, 5), but also (x, y, z)=(-3, 4, 5) is a 'virtual' primitive Pythagorean triple, holding $\rm y{=}x{+}7$ and $\rm x^2{+}y^2{=}z^2$.

Considering this duality for z=5, we put $\rm a_0{=}5$, $\rm a_1{=}17$.

Similarly as the case $\rm a_1{=}13$, substituting n=0 and n=1 into (3), we obtain a linear simultaneous equation for a and b.\\
$$ \left\{
\begin{aligned}
\rm a & \ \rm {+} & \rm b & = & \rm a_0 = 5\ \ \ \ \\
\rm a{x_{+}} & \ \rm {+} & \rm b{{x_{-}}} & = & \rm a_1 = 17\ \ 
\end{aligned}
\right.
$$
\qquad
Solving this equation, 

$\rm \ \ \ \ \ \ \ a={\dfrac{34{-}5A{+}5\sqrt{A^2{+}4B}}{2\sqrt{A^2{+}4B}}}$, $\rm b={\dfrac{5A{-}34{+}5\sqrt{A^2{+}4B}}{2\sqrt{A^2{+}4B}}}$

Then the following express holds.

$\rm {\ \ \ \ \ \ \ \ \ \ \ \ \ \ 2\ebrak{\dfrac{\rm 25{A^2{+}4B}{-}(5A{-}26)^2}{\rm {A^2{+}4B}}}\ebrak{\rm {-}B}^n{=}49}$

For keeping constancy B must be ${-}1$, and A=6, $\dfrac{386}{49}$.

If A=$\dfrac{386}{49}$ then

$\ \ \ \ \ \ \ \ \rm a{=}\dfrac{5\sqrt{2}{-}1}{2\sqrt{2}}, b{=}\dfrac{5\sqrt{2}{+}1}{2\sqrt{2}}, x_{+}{=}\dfrac{193{+}132\sqrt{2}}{49}, x_{-}{=}\dfrac{193{-}132\sqrt{2}}{49}.$\\

So $\rm a_2{=}ax_{+}^2{+}bx_{-}^2$

$\ \ \ \ \ \ \ \ \ \rm {=}\ebrak{\dfrac{5\sqrt{2}{-}1}{2\sqrt{2}}}\ebrak{\dfrac{193{+}132\sqrt{2}}{49}}^2{+}\ebrak{\dfrac{5\sqrt{2}{+}1}{2\sqrt{2}}}\ebrak{\dfrac{193{-}132\sqrt{2}}{49}}^2$

$\ \ \ \ \ \ \ \ \ \rm {=}{\dfrac{309533}{2401}}\ {\notin}{\rm N^{+}}$; A=$\dfrac{386}{49}$ is invalid.\\

If A=6 then
\begin{eqnarray}
\lefteqn {\rm \ \ \ \ a_{n{+}1}{=}6a_n{-}a_{n{-}1}.} 
\end{eqnarray}
\qquad Now evidently any $\rm a_n$ is a positive integer. And we easily see that $\rm a_2{=}6a_1{-}a_0{=}6{\times}17{-}5{=}97$, and that $97^2{=}65^2{+}72^2{=}65^2{+}(65{+}7)^2$, which results in that (65,72,97) is a primitive Pythagorean triple and an element of the set shown at (1). 97 is also suitable.

In addition, 
\begin{eqnarray}
\lefteqn {\rm  {\ {d_n}^2{=}2{a_n}^2{-}49}{=}\ebrak{\ebrak{\dfrac{5\sqrt2{+}1}{2}}(3{+}2\sqrt2)^n{+}\ebrak{\dfrac{5\sqrt2{-}1}{2}}(3{-}2\sqrt2)^n}^2{-}49}\nonumber \\ 
\lefteqn {\rm  \ \ \ \ \ \ \ {=}\ebrak{\ebrak{\dfrac{5\sqrt2{+}1}{2}}(3{+}2\sqrt2)^n{-}\ebrak{\dfrac{5\sqrt2{-}1}{2}}(3{-}2\sqrt2)^n}^2}
\end{eqnarray}

If n${\geq}$1 then $\rm{d_n}^2$ is the square of a positive integer greater than 7.\\

From the above these two recurrence relations with different initial values constitute a set of positive integers. 

Now we check if this set is identical with the set consisting of z's of the set of primitive Pythagorean triples at (1).
\\\\
2. Reaeach by using parameters.\\

We consider a generic triple (m,$\rm m{+}7$,c), instead of (m,$\rm m{+}7$,$\rm a_n$).

According to number theory, there exist r, s$\rm {\in}N^{+}$ relatively prime with different parity such that c=$\rm r^2{+}s^2$ and r>s. If and only if r and s has a common divisor greater than 1, then $\rm r^2{-}s^2$, $\rm 2rs$, $\rm r^2{+}s^2$ also have a common divisor greater than 1, which contradicts the assumption.

There are two cases. One is that m is even, then m=2rs and m+7=$\rm r^2{-}s^2{=}2rs{+}7$. The other is that m is odd, then m=$\rm r^2{-}s^2{=}2rs{+}7$ and m+7=2rs=$\rm r^2{-}s^2{+}7$.

We consider an equation
\begin{eqnarray}
\lefteqn {\rm {\ \ \ \ \ \ {r_k^2{-}s_k^2{-}2r_ks_k}{=}7{(-1)}^k}} 
\end{eqnarray}
and search a sequence of $\rm \{ (r_k,s_k) \}_{n=1,2,\cdots{,}\infty}$.

This equation covers the above-mentioned two cases and can be transformed into
\begin{eqnarray}
\lefteqn {\rm {\ \ \ \ \ \ {(r_k{-}s_k)^2{-}2{s_k}^2}{=}7{(-1)}^k}} \nonumber
\end{eqnarray}
\qquad Now we set $\rm p_k{=}r_k{-}s_k,\ q_k{=}s_k$.

We have the equation
\begin{eqnarray}
\lefteqn {\rm {\ \ \ \ \ \ {{p_k}^2{-}2{q_k}^2}{=}7{(-1)}^k}.}
\end{eqnarray}

Now we consider two above-mentioned sequences, $\rm \{13,73,\cdots \}$ and $\rm \{17,97,\cdots \}$ and apply for (11).\\

First we search $\rm \{13,73,\cdots \}$.

For z=13 at (1) the Pythagorean triple is (5,12,13), then ($\rm r_1$,$\rm s_1$)=(3,2).

And for z=73 the Pythagorean triple is (48,55,73), then ($\rm r_2$,$\rm s_2$)=(8,3).

By transforming, we obtain ($\rm p_1$,$\rm q_1$)=(1,2) and ($\rm p_2$,$\rm q_2$)=(5,3).

First we search $\rm \{17,97,\cdots \}$.

For z=17 at (1) the Pythagorean triple is (8,15,17), then ($\rm r_1$,$\rm s_1$)=(4,1), and for z=97 the Pythogorean triple is (65,72,97), then ($\rm r_2$,$\rm s_2$)=(9,4).

By transforming, we obtain ($\rm p_1$,$\rm q_1$)=(3,1) and ($\rm p_2$,$\rm q_2$)=(5,4).\\

For the former, if we define a function
\begin{eqnarray}
\lefteqn {\rm {\ \ \ \ {f\ {:}\ ({p_i},{q_i}){\longrightarrow}{p_i}{+}\sqrt2{q_i}}}}
\end{eqnarray}
\qquad Then f(${\rm p_1}$,${\rm q_1}$)=f(1,2)=1+2$\sqrt2$ and f(${\rm p_2}$,${\rm q_2}$)=f(5,3)=5+3$\sqrt2$.

5+3$\sqrt2$=(1+$\sqrt2$)$\cdot$(1+2$\sqrt2$), based on which we assume ${\rm p_1}$ and ${\rm q_1}$ such that 
\begin{eqnarray}
\lefteqn {\rm \ \ f({p_k},{q_k}){=}{p_k}{+}\sqrt2{q_k}} \nonumber \\
\lefteqn {\rm {\ \ \ \ \ \ \ \ \ \ \ \ \ \ {{=}{(1{+}\sqrt2)}^{k{-}1}f({p_1},{q_1}){=}{(1{+}\sqrt2)}^{k{-}1}(1{+}2\sqrt2)}}} \\
\lefteqn {\rm {\ \ \ \ \ \ \ \ \ \ \ \ \ \ {{=}{(1{+}\sqrt2)}^k{(3{-}\sqrt2)}}}} \nonumber
\end{eqnarray}

We can set the initial value as $\rm ({p_0},{q_0}){=}(3,-1)$.\\

Similarly for the latter, by using (12),  f(${\rm p_1}$,${\rm q_1}$)=f(3,1)=3+$\sqrt2$ and f(${\rm p_2}$,${\rm q_2}$)=f(5,4)=5+4$\sqrt2$.\

5+4$\sqrt2$=(1+$\sqrt2$)$\cdot$(3+$\sqrt2$), alike (13) we assume as follows.
\begin{eqnarray}
\lefteqn {\rm \ \ f({p_k},{q_k}){=}{p_k}{+}\sqrt2{q_k}} \nonumber \\
\lefteqn {\rm {\ \ \ \ \ \ \ \ \ \ \ \ \ \ {{=}{(1{+}\sqrt2)}^{k{-}1}f({p_1},{q_1}){=}{(1{+}\sqrt2)}^{k{-}1}(3{+}\sqrt2)}}} \\
\lefteqn {\rm {\ \ \ \ \ \ \ \ \ \ \ \ \ \ {{=}{(1{+}\sqrt2)}^k{({-}1{+}2\sqrt2)}}}} \nonumber
\end{eqnarray}

We can set the initial value as $\rm ({p_0},{q_0}){=}(-1,2)$.\\

Without loss of generality arranging (11) into 2 equations, 
\begin{eqnarray}
\lefteqn {\rm {\ \ \ \ \ \ {{p_k}^2{-}2{q_k}^2}{=}7{(-1)}^k}.} \\
\lefteqn {\rm {\ \ \ \ \ \ {{p_k}^2{-}2{q_k}^2}{=}{-}7{(-1)}^k}.} 
\end{eqnarray}

Any ($\rm{p_k}$,$\rm{q_k}$) of (12) fulfills (16), and that of (13) fulfills (15).\\

More generalizing, this problem for now is to find any solution $\rm(p,q)$ for
\begin{eqnarray}
\lefteqn {\rm {\ \ \ \ \ \ {p^2{-}2q^2}{=}{\pm}7,\ where\ p\ and\ q\ are\ positive \ integers.}} 
\end{eqnarray}

Now assume (x,y) fulfills (17).

Then ${\rm x^2{-}2y^2}{=}{\pm}7$.

${\rm \ \ (x{+}\sqrt2{y})({1}{+}\sqrt2){=}(x{+}2y){+}(x{+}y)\sqrt2}$, ${\rm \dfrac{{x}{+}\sqrt2{y}}{{1}{+}\sqrt2}{=}({-}x{+}2y){+}(x{-}y)\sqrt2}$

So ${\rm (x{+}2y)^2{-}2{(x{+}y)}^2{=}{-}x^2{+}2y^2{=}{\mp}7}$, and ${\rm ({-}x{+}2y)^2{-}2{(x{-}y)}^2{=}{-}x^2{+}2y^2{=}{\mp}7}$.

Moreover, we define ($\rm x_s$,$\rm y_s$) and ($\rm x_t$,$\rm y_t$) as
\begin{eqnarray}
\lefteqn {\rm x_s{+}\sqrt2{y_s}{=}{(1{+}\sqrt2)}^{s}(x{+}\sqrt2{y})}\\
\lefteqn {\rm x_t{+}\sqrt2{y_t}{=}\dfrac{(x{+}\sqrt2{y})}{(1{+}\sqrt2)^{t}}}\\
\lefteqn {\rm \ \ where\ s\ and\ t\ are\ positive\ integers.} \nonumber 
\end{eqnarray}

Then inductively for any s and t, ($\rm x_s$,$\rm y_s$) and ($\rm x_t$,$\rm y_t$) fulfills the equation (17), \\
under the expanded condition that p and q are integers (p and q need not be positive) and ${\rm p{+}\sqrt2{q}>0}$.

Therefore any solution of (17) has its derivative solution in the form of  ($\rm x_s$,$\rm y_s$) or  ($\rm x_t$,$\rm y_t$) or itself, having a value between ${1}{+}\sqrt2$ and $({1}{+}\sqrt2)^2$=${3}{+}2\sqrt2$.

Now we will find locally the above-mentioned solutions, between ${1}{+}\sqrt2$ and ${3}{+}2\sqrt2$. Note that p is odd. 
\\\\
<< For the case p,q>0 >>

There are 5 numbers eligible for ${\rm p{+}\sqrt2{q}}$ whose values are between ${1}{+}\sqrt2$ and ${3}{+}2\sqrt2$. They are as follows, ${1}{+}\sqrt2$, ${3}{+}\sqrt2$, ${1}{+}2\sqrt2$, ${3}{+}2\sqrt2$ and ${3}{+}2\sqrt2$.

With simple computations we see that ${1}{+}2\sqrt2$ and ${3}{+}\sqrt2$ hold the equation (17).
\\\\
<< For the case p=0 or q=0 >>

${\rm {p^2{-}2q^2}{\neq}{\pm}7}$, so invalid.
\\\\
<< For the case p, q<0 >>

${\rm p{+}\sqrt2{q}<0}$, so invalid.
\\\\
<< For the case p>0 or q<0 >>

${\rm p{+}\sqrt2{q}{=}\sqrt{2q^2\pm7}{+}\sqrt2{q}{=}\sqrt{2{\left| q \right|}^2\pm7}{-}\sqrt2{\left| q \right|}}$\\
${\rm Now\ p{+}\sqrt2{q}>0,\ so}$\\
${\rm \ \ \ \ \ \ p{+}\sqrt2{q}{=}\sqrt{2{\left| q \right|}^2{+}7}{-}\sqrt2{\left| q \right|}{=}\dfrac{7}{\sqrt{2{\left| q \right|}^2{+}7}{+}\sqrt2{\left| q \right|}}{\le}\dfrac{7}{3{+}\sqrt2}{=}3{-}\sqrt2<1{+}\sqrt2}$, so invalid.
\\\\
<< For the case p<0 or q>0 >>

${\rm p{+}\sqrt2{q}{=}p{+}\sqrt{p^2\pm7}{=}\sqrt{{\left| p \right|}^2\pm7}{-}{\left| p \right|}}$\\
${\rm Now\ p{+}\sqrt2{q}>0,\ so}$\\
${\rm \ \ \ \ \ \ p{+}\sqrt2{q}{=}\sqrt{2{\left| p \right|}^2{+}7}{-}\sqrt2{\left| p \right|}{=}\dfrac{7}{\sqrt{2{\left| p \right|}^2{+}7}{+}\sqrt2{\left| p \right|}}{\le}\dfrac{7}{3{+}\sqrt2}{=}3{-}\sqrt2<1{+}\sqrt2}$, so invalid.
\\\\

Overall only ${1}{+}2\sqrt2$ and ${3}{+}\sqrt2$ suit the equation (17) between ${1}{+}\sqrt2$ and ${3}{+}2\sqrt2$. 

In other words, nothing but (1,2) and (3,1) are the solution of (17) under the restriction of number size above.

Therefore (13) and (14) cover all solutions of (1).

Developing (13), 

$\rm \ \ \ \ {p_k}{+}\sqrt2{q_k}{=}{(3{-}\sqrt2)}{(1{+}\sqrt2)}^k$

$\rm \ \ \ \ {p_k}{-}\sqrt2{q_k}{=}{(3{+}\sqrt2)}{(1{-}\sqrt2)}^k$
\begin{eqnarray}
\lefteqn{ \rm \ \ \ \ {q_k}{=}{s_k}{=}{\dfrac{1}{4}}\ebrak{\rm {\sqrt2(3{-}\sqrt2)}{(1{+}\sqrt2)}^{k}{-}{\sqrt2(3{+}\sqrt2)}{(1{-}\sqrt2)}^{k}}}\nonumber \\
\lefteqn{ \rm \ \ \ \ {p_k}{=}{\dfrac{1}{2}}\ebrak{\rm {(3{-}\sqrt2)}{(1{+}\sqrt2)}^{k}{+}{(3{+}\sqrt2)}{(1{-}\sqrt2)}^{k}}}\\
\lefteqn{ \rm \ \ \ \ {r_k}{=}{\dfrac{1}{4}}\ebrak{\rm {(4{+}\sqrt2)}{(1{+}\sqrt2)}^{k}{+}{(4{-}\sqrt2)}{(1{-}\sqrt2)}^{k}}} \nonumber
\end{eqnarray}

$\rm {s_k}{=}{\dfrac{1}{4}}\ebrak{\rm {\sqrt2(3{-}\sqrt2)}{(1{+}\sqrt2)}^{k}{-}{\sqrt2(3{+}\sqrt2)}{(1{-}\sqrt2)}^{k}}$

$\rm \ \ \ \ {=}{\dfrac{1}{4}}\ebrak{\rm {\sqrt2(3{-}\sqrt2)}{(1{+}\sqrt2)}{(1{+}\sqrt2)}^{k{-}1}{-}{\sqrt2(3{+}\sqrt2)}{(1{-}\sqrt2)}{(1{-}\sqrt2)}^{k{-}1}}$

$\rm \ \ \ \ {=}{\dfrac{1}{4}}\ebrak{\rm {(4{+}\sqrt2)}{(1{+}\sqrt2)}^{k{-}1}{+}{(4{-}\sqrt2)}{(1{-}\sqrt2)}^{k{-}1}}$

$\rm \ \ \ \ {=}{r_{k{-}1}}$.
\\

The equations above shows that $\rm {r_k}{=}2{r_{k{-}1}}{+}{r_{k{-}2}}$ and $\rm {s_k}{=}2{s_{k{-}1}}{+}{s_{k{-}2}}$.

So

$\rm \ \ \ \ \ \ {r_k}^2{+}{s_k}^2{=}{s_{k{+}1}^2{+}{s_k}^2{=}(2{s_k}{+}{s_{k{-}1}})^2{+}{s_k}^2{=}5{s_k}^2{+}4{s_k}{s_{k{-}1}}{+}{s_{k{-}1}}^2  }$

$\rm \ \ \ \ \ \ \ \ \ \ \ \ \ \ \ \ \ \ \ {{=}6{s_k}^2{+}5{s_{k{-}1}}^2{-}{(2s_{k{-}1}{-}s_k)}^2{=}6{s_k}^2{+}5{s_{k{-}1}}^2{-}{s_{k{-}2}}^2}$

$\rm \ \ \ \ \ \ \ \ \ \ \ \ \ \ \ \ \ \ \ {{=}6({s_k}^2{+}{s_{k{-}1}}^2){-}({s_{k{-}1}}^2{+}{s_{k{-}2}}^2)  }$

$\rm \ \ \ \ \ \ \ \ \ \ \ \ \ \ \ \ \ \ \ {{=}6({r_{k{-}1}}^2{+}{s_{k{-}1}}^2){-}({r_{k{-}2}}^2{+}{s_{k{-}2}}^2)  }$ .\\

$\rm \therefore {c_k}{=}6{c_{k{-}1}}{-}{c_{k{-}2}}$\\

For (6) and (13), $\rm {a_0}{=}{c_0}{=}5$ and $\rm {a_1}{=}{c_1}{=}13$.

So we see that the two sets, $\rm \{a_n\}_{n=1,2,\cdots{,}\infty}$ for (6) and $\rm \{c_n\}_{n=1,2,\cdots{,}\infty}$ for now, are term by term identical.\\

Similarly for (14)
\begin{eqnarray}
\lefteqn{ \rm \ \ \ \ {q_k}{=}{s_k}{=}{\dfrac{1}{4}}\ebrak{\rm {\sqrt2(-1{+}2\sqrt2)}{(1{+}\sqrt2)}^{k}{-}{\sqrt2(-1{-}2\sqrt2)}{(1{-}\sqrt2)}^{k}}}\nonumber \\
\lefteqn{ \rm \ \ \ \ {p_k}{=}{\dfrac{1}{2}}\ebrak{\rm {(-1{+}2\sqrt2)}{(1{+}\sqrt2)}^{k}{+}{(-1{-}2\sqrt2)}{(1{-}\sqrt2)}^{k}}}\\
\lefteqn{ \rm \ \ \ \ {r_k}{=}{\dfrac{1}{4}}\ebrak{\rm {(2{+}3\sqrt2)}{(1{+}\sqrt2)}^{k}{+}{(2{-}3\sqrt2)}{(1{-}\sqrt2)}^{k}}} \nonumber
\end{eqnarray}

And then equations, $\rm {r_k}{=}{s_k}$ and $\rm {c_k}{=}6{c_{k{-}1}}{-}{c_{k{-}2}}$ also hold.

For (8) and (14), $\rm {a_0}{=}{c_0}{=}5$ and $\rm {a_1}{=}{c_1}{=}17$.

Again we see that the two sets, $\rm \{a_n\}_{n=1,2,\cdots{,}\infty}$ for (8) and $\rm \{c_n\}_{n=1,2,\cdots{,}\infty}$ for now, are also term by term identical.\\

As a result the union of two sets, 
\begin{eqnarray}
\lefteqn{ \rm \{a_n\}_{n=1,2,\cdots{,}\infty}\ where\ {a_{n{+}1}}{=}6{a_n}{-}{a_{n{-}1}}\ with\ {a_0}{=}5\ and\ {a_1}{=}13} \\
\lefteqn{\rm \ \ \ \ \ \ and} \nonumber \\
\lefteqn{ \rm \{a_n\}_{n=1,2,\cdots{,}\infty}\ where\ {a_{n{+}1}}{=}6{a_n}{-}{a_{n{-}1}}\ with\ {a_0}{=}5\ and\ {a_1}{=}17} 
\end{eqnarray}
coincides with the set of z's at (1), where any z is the greatest number of an primary Pythagorean triple.
\\\\
\noindent 3. Discussion and Conclusions.\\

As mentioned above we have showed that the set, consisting of the greatest number of any element, which is a primitive Pythagorean triple, of the set described at (2), coincides with two sequences (22) and (23) each composed of a similar recurrence relation.

Cimmino [1] treated the set of maximum numbers of all primitive Pythagorean triples such that the difference between the two non-maximum numbers is 1. 

And we treat now for the case the difference is 7.

When the difference is 1, the set of z's can be written with only one recurrence relation. But for our case, the difference is 7 and the set needs as many as two recurrence relations, (22) and (23).

The transitions of the terms in the two sequences (22) and (23) can be described as follows.

\begin{figure}[htbp]
\begin{center}
\includegraphics [width=80mm]{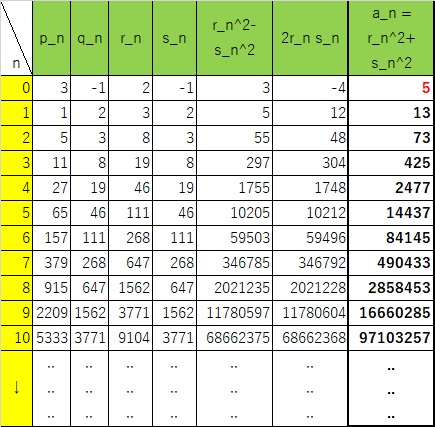}
\caption{For the case $\rm a_1$=13}
\end{center}
\end{figure}
\begin{figure}[htbp]
\begin{center}
\includegraphics [width=80mm]{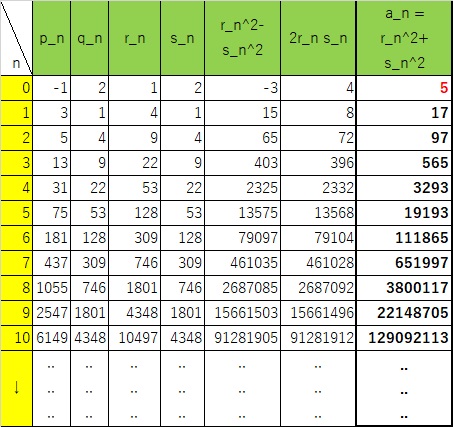}
\caption{For the case $\rm a_1$=17}
\end{center}
\end{figure}
\newpage
Though controlled by one of the two recurrence relations (22) and (23) which are almost similar to each other, the two sequences look completely different from each other, as we see as Figure 1 and Figure 2.

At a glance the two sequences seems to have an entirely different stream from each other. \\

Now we expand the range of n to negative.\\
\begin{eqnarray}
\lefteqn {\rm \ \ {p_{k{-}1}}{+}\sqrt2{q_{k{-}1}}{=}  f({p_{k{-}1},{q_{k{-}1})}}} \nonumber \\
\lefteqn {\rm {\ \ \ \ \ \ \ \ \ \ \ \ \ \ {{=}{(1{+}\sqrt2)}^{{-}1}f({p_k},{q_k}){=}{({-}1{+}\sqrt2)}({p_k}{+}{q_k}\sqrt2)}}} \nonumber \\
\lefteqn {\rm {\ \ \ \ \ \ \ \ \ \ \ \ \ \ {{=}{(2{q_k}{-}{p_k}){-}\sqrt2({p_k}{-}{q_k})}}}}\nonumber \\
\lefteqn {\rm {\ \ \ \ \ \ \ \ \ \ ,\ so\ evidently}} \nonumber \\
\lefteqn {\rm {\ \ \ \ \ \ \ \ \ \ \ \ \ \ {p_{k{-}1}}{=}2{q_k}{-}{p_k},\ {q_{k{-}1}}{=}{p_k}{-}{q_k} \ . }} 
\end{eqnarray}

By using (24) we obtain an expanded table for the case of (22) as below.
\begin{figure}[htbp]
\begin{center}
\includegraphics [width=80mm]{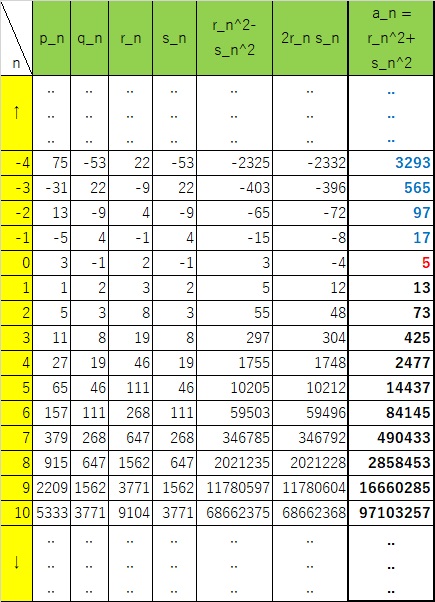}
\caption{Expansion for the case $\rm a_1$=13}
\end{center}
\end{figure}
\newpage
According to Figure 3, $\rm a_{{-}1}{=}17$, $\rm a_{{-}2}{=}97$, $\rm a_{{-}3}{=}565$ and $\rm a_{{-}4}{=}3293$.

And on Figure 2, which shows the sequence $\rm \{a_n\}_{n=1,2,\cdots{,}\infty}$ for the case $\rm a_1{=}17$, we see $\rm a_{1}{=}17$, $\rm a_{2}{=}97$, $\rm a_{3}{=}565$ and $\rm a_{4}{=}3293$. 

We can assume that $\rm a_{{-}n}$ for (22) is equal to $\rm a_{n}$ for (23).

By the way, a recurrence relation $\rm a_{n{+}1}{=}6a_n{-}a_{n{-}1}$ is transformed into 
\begin{eqnarray}
\lefteqn {\rm \ \ a_{n{-}1}{=}6a_n{-}a_{n{+}1}} 
\end{eqnarray}

Another recurrence relation (25) of the same form has appeared, so $\rm a_{{-}n}$ for (22) is equal to $\rm a_{n}$ for (23), for any n>0. 

Thus the above-mentioned assumption has inductively proved to be true.\\

Similarly by using (24) we get another expanded table for the case of (23).
\begin{figure}[htbp]
\begin{center}
\includegraphics [width=80mm]{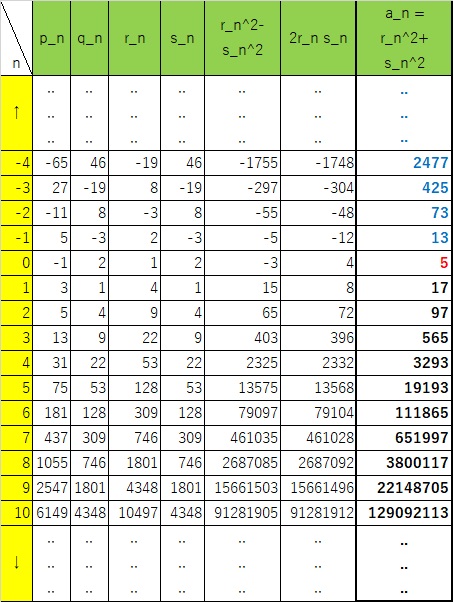}
\caption{Expansion for the case $\rm a_1$=17}
\end{center}
\end{figure}
\newpage
According to Figure 4, $\rm a_{{-}1}{=}13$, $\rm a_{{-}2}{=}73$, $\rm a_{{-}3}{=}425$ and $\rm a_{{-}4}{=}2477$.

And on Figure 1, which shows the sequence $\rm \{a_n\}_{n=1,2,\cdots{,}\infty}$ for the case $\rm a_1{=}13$, we see $\rm a_{1}{=}13$, $\rm a_{2}{=}73$, $\rm a_{3}{=}425$ and $\rm a_{4}{=}2477$.

Inductively we see that $\rm a_{{-}n}$ for (23) is equal to $\rm a_{n}$ for (22).

As a result, all z's of (1) can be enumerated by a sequence composed of a single recurrence relation; both Figure 3 and Figure 4 include such a good sequence that its latest value changes in a reverse direction to each other.

In the connection, we show expanded table for the case of Cimmino's [1].
\begin{figure}[htbp]
\begin{center}
\includegraphics [width=80mm]{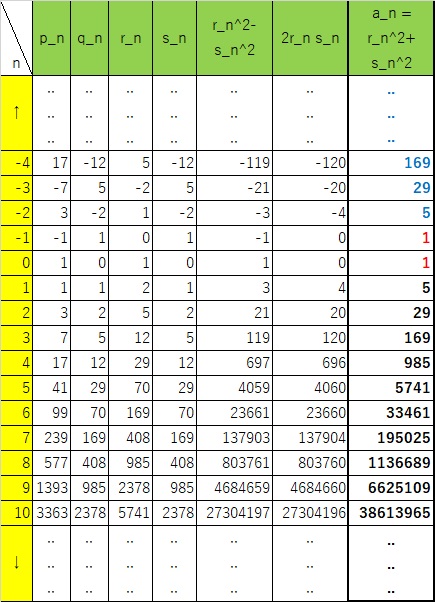}
\caption{Expansion for the case of Cimmino's}
\end{center}
\end{figure}
\newpage
\centerline{References}
\quad\par
\noindent [1] Luigi Cimmino, arxiv, available at \href{https://arxiv.org/abs/math/0510417}{https://arxiv.org/abs/math/0510417}\\
\end{document}